\documentclass[letterpaper, 10 pt, conference]{ieeeconf}

\IEEEoverridecommandlockouts    
\usepackage{cite}
\usepackage{amsmath,amssymb,amsfonts,lipsum}
\usepackage{algorithmic}
\usepackage{graphicx}
\usepackage{textcomp}
\usepackage{comment}
\usepackage{color}

\usepackage{graphicx}
\newtheorem{ps}{Problem}
\newtheorem{asump}{Assumption}
\newtheorem{theorem}{Theorem}
\newtheorem{remark}{Remark}

\newtheorem{lemma}{Lemma}

\title{Analog Data-Driven Theory and Estimation of the Region of Attraction Using Sampled-Data}
\author{Karthik Shenoy$^{1}$, Arvind Ragghav$^{2}$, Vijaysekhar Chellaboina$^{3}$
\thanks{$^{1}$KS is a Prime Minister's Research Fellowship scholar at the Dept. of Electrical Engineering, IIT-Madras, India.
        {\tt\small ee21d405@smail.iitm.ac.in}}
\thanks{$^{2}$ AR is a Project Associate at ICS\&R, IIT Madras, India.
        {\tt\small arvind00ragg@gmail.com}}
\thanks{$^{3}$ VC is the Dean, School of Computer Science, UPES Dehradun. 
        {\tt\small vijaysekhar@alumni.iitm.ac.in }}%
}

\begin{document}

\maketitle
\thispagestyle{empty}
\pagestyle{empty}
\begin{abstract}
The contributions of this technical note are twofold. Firstly, we formulate an optimization problem to obtain a linear representation of a nonlinear vector field based on a system's trajectory. We also prove that its cost function is strictly convex, given the trajectory is persistently exciting. Under certain observability conditions, we provide results that guarantee the Hurwitz stability of the global minimizer. Secondly, we present a novel algorithm based on point-wise geometric flows to estimate the boundary of the region of attraction. We show that the algorithm converges to the exact boundary of the region of attraction under certain assumptions on the system dynamics. Finally, we validate the results using simulations on various nonlinear autonomous systems. 
\end{abstract}

\section{Introduction}
\label{sec:introduction}
Data-driven control and optimization has been one of the most extensively researched areas in modern control theory during the past decade. The lack of reliable models for complex systems and the availability of vast amounts of process data have led to this surge in using data-driven techniques to design controllers. 
These techniques, based on the works by Willems et al.,\cite{POE}, directly use \textit{persistently exciting} input-output data obtained from a plant to design control laws, eliminating the need for an intermediate system identification stage. While the preponderance of literature is dedicated to data-driven controller design for linear dynamical systems \cite{formulabased, informativity,mpc, minen, delay, dissipativity}, remarkable contributions have been made in the areas of data-driven control design for networks \cite{pasqdatadrivennetworks},\cite{pasqnetworknew}, and nonlinear systems \cite{NLFBdiction, NLdatadriven,tabuadafbdata} to name a few. One can also get an overview of all the aforementioned articles in the recent survey paper \cite{DDsurvey}. In all the mentioned articles, the analysis is carried out mainly in the discrete-time setting or using sampled data for continuous-time systems. Moreover, these analyses are often restricted to a local region around the equilibrium point of a nonlinear system. To our knowledge, data-driven analysis using analog trajectory data for continuous-time nonlinear systems, without sampling, is scant. 

On another note, estimating the region of attraction (ROA) of an equilibrium point, of a nonlinear autonomous system, has received widespread attention since the works by Genesio et al., \cite{roa_seminal}. Utilizing Lyapunov functions to estimate the ROA has been discussed for nonlinear systems in \cite{roa_maximal_lyapunov} and robust estimation of the ROA in \cite{roa_parameter}. In \cite{roa_estimate_LMI}, the authors formulate this estimation problem as a solution to a Linear Matrix Inequality for uncertain polynomial systems. Trajectory-based methods have been employed in \cite{roa_trajectory} to estimate an intermediate Lyapunov function and in \cite{roa_trajectory_forward_reachable} by computing the forward reachable sets. Recently, data-driven methods for estimating the ROA have gained popularity. The authors in \cite{roa_operator} and \cite{roa_neural}, estimate the ROA using operator-theoretic methods and neural-network-based methods respectively. Most of the methods, as mentioned earlier, rely on constructing a Lyapunov function as an intermediate step or using trajectory data to construct forward reachable sets to estimate the ROA. Furthermore, several of the previously mentioned approaches can only provide a conservative estimate of the ROA.

Thus, motivated by both these premises, we analyze data-driven methods for nonlinear systems in the analog setting in this article without resorting to sampled data and when the trajectories are not localized around the origin. Further, we address the problem of obtaining the ROA from data when the system model is unavailable. 

\textit{Contributions}. We encapsulate the main contributions of this article below
\begin{itemize}


    \item We propose an optimization problem to obtain a linear vector field that best matches the vector field of an autonomous dynamical system along a given trajectory. The objective function in the optimization problem is an integral cost function, which is dependent on the trajectories of the system under consideration. We prove that the cost function is strictly convex if the trajectory is persistently exciting  (Theorem \ref{thm:convex}). Hence, we propose a gradient-flow dynamics, to solve the problem.
    
    \item We prove that the global minimizer of the cost function is Lyapunov stable if the initial condition of the system trajectory belongs to the region of attraction (Theorem \ref{thm:stability}). Furthermore, under certain observability conditions, the global minimizer is shown to be Hurwitz.
    \item We propose a novel algorithm, based on point-wise geometric flows, to estimate the region of attraction of a locally exponentially stable equilibrium point for a general class of nonlinear autonomous systems, using the minimizer of the optimization problem. We also discuss the algorithm's convergence and prove that it converges to the exact boundary of the ROA for any given initial conservative guess of the ROA under certain radial unboundedness conditions on the dynamics.   
\end{itemize}

The rest of the letter is organized as follows. We define the optimization problem in Section \ref{sec: ProblemDef}. Section \ref{sec:res_analysis} discusses its solution and analysis. Section \ref{sec:results_ROA} is dedicated to estimating the ROA. Finally, in Section \ref{sec: simulations}, we validate the proposed algorithm through simulations on various systems whose ROA exhibits different geometric properties.

\textit{Notation}. We denote the space of $n$-dimensional real vectors by $\mathbb{R}^n$ and the $i^{\text{th}}$ component of $x\in\mathbb{R}^{n}$ as $x_{i}$. The set of whole numbers, non-negative real numbers, the extended real line, and $\{1,\hdots,n\}$ are denoted by $\mathbb{W}, \mathbb{R}_{\geq 0}, \mathbb{R}\cup\{\infty\}, [n]$, respectively. The partial derivative of a function $f$, with respect to one of its entries $x$ (which might be scalar or vector-valued), is denoted by $f_{x}$. The $2$-norm of a vector in $\mathbb{R}^n$ is denoted by $\|\cdot\|$. The Kronecker product between $A$ and $B$ is denoted by  $A \otimes B$. Given a matrix $A\in\mathbb{R}^{n\times n}$, $A\succ 0, A\preceq 0, \|A\|_{F}, A^{\dagger}$ denote that A is positive definite, negative-semi definite, its Frobenius norm, and Moore-Penrose inverse, respectively. Given a set $N$, its interior, boundary, closure, and complement are denoted by $\operatorname{int}(N),\partial N,\;\overline{N},\;N^c$, respectively. We write $f(x)=o(g(x))$ to denote that$\|f(x)\|/\|g(x)\|\to0$ as $\|x\|\to 0$. 
\section{Problem Definition}
\label{sec: ProblemDef}
Consider the nonlinear autonomous system
\begin{align}
    \dot{x}(t)=f(x(t)),\quad x(0)=x_0\label{eqn:autosystem}
\end{align}
where $f:\;\mathbb{R}^n\to\mathbb{R}^n$ is smooth. Without any loss of generality, let $x^*=0$ be a locally exponentially stable equilibrium point of \eqref{eqn:autosystem}, with a region of attraction $M\subseteq \mathbb{R}^n$. Let $A\triangleq f_x(0)$ be the Jacobian matrix computed at the origin. The matrix $A\in\mathbb{R}^{n\times n}$ is Hurwitz as a consequence of the origin being locally exponentially stable. Now we re-write the dynamics \eqref{eqn:autosystem} as 
\begin{align}
     \dot{x}(t)=Ax(t)+g(x(t)),\;x(0)=x_0\label{eqn:autosystem_lin_ho}
\end{align}
where $g(x)$ defined as $f(x)-Ax$ represents the higher order terms in $f(x)$. We assume that $g(x)=o(x)$, which implies that $\exists\; D_0\subset \mathbb{R}^n$ containing the origin and a constant $\gamma>0$ such that $\forall\;x\in D_0,\;\|g(x)\|\leq\gamma\|x\|$. Let $s(t,x_0)$ be the solution of \eqref{eqn:autosystem}.  Then we define a cost $J:\;\mathbb{R}^{n\times n}\times \mathbb{R}^n\to \mathbb{R}_{\geq 0}\cup\{\infty\}$ as
\begin{align}
    J(\bar{A},\;x_0)=\int_0^\infty\|f(s(t,x_0))-\bar{A}s(t,x_0)\|^2 \mathrm{d}t\label{eqn:cost}
\end{align}
associated with a given matrix $\bar{A}$ and an initial condition $x_0\in M$. We assume the following in this article.
\begin{asump}[Persistency of Excitation]\label{assm:POE}
    The matrix $\int_0^\tau s(t,x_0)s(t,x_0)^\top\mathrm{d}t$ is positive definite $\forall\;\tau\in\mathbb{R}_{\geq 0}$ and $\forall\,x_0\in\mathbb{R}^n$.
\end{asump}

Based on this assumption, we address the following problems in this article.

\begin{ps}\label{prob:1}
    \textit{For a nonlinear system of the form \eqref{eqn:autosystem}, with $f(x)$ unknown and given a trajectory $s(t,x_0)$, estimate $\hat{A}(x_0)$ satisfying  
    \begin{align}
        \hat{A}(x_0)=\underset{\bar{A}\;\in\;\mathbb{R}^{n\times n}}{\operatorname{argmin}}\;J(\bar{A},\;x_0)
    \end{align}
    which gives the best linear fit $\hat{A}(x_0)s(t,x_0)$ of $f(s(t,x_0))$.} \label{P1}
\end{ps}

\begin{ps}\label{prob:2}
    \textit{Estimate the region of attraction $M$ of the origin, under the dynamics \eqref{eqn:autosystem}, using trajectory data.} when the function $f(x)$ is unknown.
\end{ps}
 
\section{Results: Linear Estimate}
\label{sec:res_analysis}

We first show that the cost given in \eqref{eqn:cost} is finite for a given $x_0\in M$ and for any arbitrary matrix $\bar{A}\in\mathbb{R}^{n\times n}$.
\begin{lemma}
    $J(\bar{A},x_0)<\infty$ for any given $x_0\in M$ and $\bar{A}\in\mathbb{R}^{n\times n}$.
\end{lemma}

\begin{proof}  
We omit the arguments in $s(t,x_0)$ due to space restrictions in all the proofs unless it is required. Consider the cost 
    \begin{align*}
     J(\bar{A}, x_{0})&=\int_0^\infty\|f(s)-\bar{A}s\|^2 \mathrm{d}t\\ 
     &=\int_0^\infty\|(A-\bar{A})s+g(s)\|^2 \mathrm{d}t\\
     &\leq \int_0^\infty(c\|(A-\bar{A})\|_F\|s\|+\|g(s)\|)^2\mathrm{d}t.
    \end{align*}
    Since $g(s)=o(s)$, there exists a subset $D_0\subset M$, containing the origin, and a constant $\gamma$ such that $\|g(s)\|\leq \gamma \|s\|, \forall\;s\in D_0$. Since $x_0\in M$, there exists a time $T_1>0$ such that $s(t,x_0)\in D_0,\;\forall\;t \geq T_1$. Next, since the origin is a locally exponentially stable equilibrium, $\exists \; D_1\subseteq M$ containing the origin, and constants $\alpha>0,\beta>0$ such that $\|s(t,\delta)\|<\alpha e^{-\beta t}, \forall\; \delta \in D_{1}$. Since $x_0\in M$ there exists a $T_2>0$ such that $s(t,x_0)\in D_1,\; \forall\; t\geq T_2$. Then $s(t,x_0)\in D_0\cap D_1,\;\forall\;t \geq T= \max\{T_1,T_2\}$.
    Thus, splitting the integral, 
    \begin{align*}
       J(\bar{A}, x_{0})&\leq I_{T} +\int_{T}^\infty(\|(A-\bar{A})\|\|s\|+\|g(s)\|)^2\mathrm{d}t \\
       &\leq I_{T}+\int_{T}^\infty\kappa\|s\|^2\mathrm{d}t
    \end{align*}
    where $I_{T}=\int_0^{T}(\|(A-\bar{A})\|\|s\|+\|g(s)\|)^{2}\mathrm{d}t$,\\ $\kappa= (\|A - \bar{A}\| + \gamma)^{2}$. Simplifying we get
    \begin{align*}
    J(\bar{A}, x_{0}) <  I_T +\int_T^\infty \kappa\alpha^2 e^{-2\beta t}\mathrm{d}t=I_T+\frac{\kappa \alpha^2}{2\beta}e^{-2\beta T}.
    \end{align*}
    which is finite for a given $x_0$ and $\bar{A}$.  
\end{proof}

\subsection{Properties of the Cost Function and the Minimizer}
\subsubsection{Strict Convexity of the Cost Function}  Define two matrices $\Gamma_1\triangleq\int\limits_0^\infty f(s(t,x_0))s^{\top}(t,x_0)\mathrm{d}t,\;\Gamma_2\triangleq\int\limits_0^\infty s(t,x_0)s^{\top}(t, x_0)\mathrm{d}t$.
Next, we show that the cost $J(\bar{A},x_0)$ is strictly convex in $\bar{A}$ and compute the minimizer $\hat{A}(x_0)$, under the condition that the trajectory $s(t,x_0)$ is \textit{persistently exciting} in the state-space, i.e., when assumption \ref{assm:POE} holds.  
\begin{theorem}\label{thm:convex}
 The cost $J(\bar{A},\;x_0)$, given by \eqref{eqn:cost} is strictly convex in $\bar{A}$ if and only if $\Gamma_{2}\succ 0$.  Moreover, the matrix 
    \begin{align}
        \hat{A}(x_0)\overset{\Delta}{=}\Gamma_1\Gamma_2^{-1}\label{eqn:Ahat}
    \end{align}
    is the global minimizer of \eqref{eqn:cost}.
\end{theorem}

\begin{proof}
The cost given by \eqref{eqn:cost} can be modified as follows  
\begin{align*}
    J(\bar{A},x_0)&=\int_0^\infty\|f(s)-\bar{A}s\|^2\\
    &=\int_0^\infty\mathrm{Trace}[(f(s)-\bar{A}s)(f(s)-\bar{A}s)^\top]\mathrm{d}t\\
    &=\mathrm{Trace}[\int_0^\infty(f(s)-\bar{A}s)(f(s)-\bar{A}s)^\top\mathrm{d}t].
\end{align*}
Expanding the outer product and completing the squares, 
\begin{align}
    J(\bar{A},x_0)=\mathrm{Trace}[&(\bar{A}-\Gamma_1\Gamma_2^{-1})\Gamma_2(\bar{A}-\Gamma_1\Gamma_2^{-1})^\top\nonumber\\&+\Gamma_0-\Gamma_1\Gamma_2^{-1}\Gamma_1^\top)]\label{eqn:cost_modified}
\end{align}
where $\Gamma_0=\int_0^\infty f(s)f(s)^\top\mathrm{d}t$.\\ 
\textbf{Sufficiency Condition:} Let $\Gamma_2\succ 0$. We only have to show that the function $\mathrm{Trace}[(\bar{A}-\Gamma_1\Gamma_2^{-1})\Gamma_2(\bar{A}-\Gamma_1\Gamma_2^{-1})^\top]$ is strictly convex in $\bar{A}$. This can be shown by invoking \cite[Fact 8.14.16]{bernstein_matrix}. Next, we prove the necessary condition by reductio ad absurdum.\\
\textbf{Necessary Condition:} We have to show that the matrix $\Gamma_2$ does not have a zero eigenvalue, as it is at least positive semi-definite by definition. Assume the contrary that $\Gamma_2$ has a $0$ eigenvalue. Thus, $\exists\;v\in \mathbb{R}^n:\;\Gamma_2v=0$. Now define $V\in \mathbb{R}^{n\times n}$ as the matrix with $n$-copies of the vector $v$ in its columns. Then $\mathrm{Trace}[(\bar{A}-\Gamma_1\Gamma_2^{-1})\Gamma_2(\bar{A}-\Gamma_1\Gamma_2^{-1})^\top]$ at $\bar{A}=\hat{A}+tV$ is $0,\;\forall\;t\in \mathbb{R}$. Hence, the minimum value of the function $\mathrm{Trace}[(\bar{A}-\Gamma_1\Gamma_2^{-1})\Gamma_2(\bar{A}-\Gamma_1\Gamma_2^{-1})^\top]$ is attained at all points $\hat{A}+tV$. This contradicts the hypothesis that the function \eqref{eqn:cost} is strictly convex. Thus, if the cost function \eqref{eqn:cost} is strictly convex,  then $\Gamma_2$ is positive definite. Hence we have derived both the necessary and the sufficiency conditions for the strict convexity of the cost function \eqref{eqn:cost}. Moreover, one can verify that the cost function attains its minima only at $\bar{A}=\Gamma_1\Gamma_2^{-1}$ from \eqref{eqn:cost_modified}, concluding the proof. 
\end{proof}

Interested readers are encouraged to derive these results using variational techniques on the space of $n\times n$ matrices. Thus, Theorem \ref{thm:convex} enables us to obtain a system representation based on its trajectories. In fact, various articles \cite{formulabased}, \cite{informativity} on data-driven analysis and control use special cases of Theorem \ref{thm:convex} to obtain a closed-loop system representation for the design of controllers. 
\begin{remark}
    Note that $\Gamma_2\succ0$ is the same persistency of excitation condition ubiquitous in data-driven systems theory (refer \cite{POE}). This condition ensures that the trajectories sufficiently rotate in the state space. Even though most classes of autonomous nonlinear systems satisfy this condition, there can be exceptional cases where the trajectory is restricted to an invariant subspace of $\mathbb{R}^n$ (for example, any invariant subspace of a nonlinear system) for all time. In the rest of the manuscript, we only consider systems whose trajectories are persistently exciting to perform our computations.
\end{remark}

\begin{remark}[Gradient Descent in $\mathbb{R}^{n\times n}$]
    Computing the inverse of $\Gamma_2$ can be computationally expensive. Hence, we can perform gradient descent in $\mathbb{R}^{n\times n}$
    \begin{align}\label{eqn: graddesc}
        \dot{B}(\tau)=\Gamma_1-B(\tau)\Gamma_2
    \end{align}
    to compute the minimizer iteratively. This is a matrix differential equation linear in $B(\tau)$. Performing a change of coordinates given by $\bar{B}(\tau)\triangleq B(\tau)-\Gamma_1\Gamma_2^{-1}$ reduces \eqref{eqn: graddesc} to $\dot{\bar{B}}(\tau)=-\bar{B}(\tau)\Gamma_2$.
    Since $\Gamma_2\succ0$,  $-\Gamma_2$ is Hurwitz. Thus $\bar{B}=0$ is a globally asymptotically stable (GAS) equilibrium point. Therefore, $B\equiv\Gamma_1\Gamma_2^{-1}$ is a GAS solution of \eqref{eqn: graddesc}.   
\end{remark}

\begin{remark}[Linear Systems]\label{rem:linearsys}
    For a linear system $\dot{x}(t)=Ax(t),\;x(0)=x_0$, where $A$ is Hurwitz, $\hat{A}(x_0)=A,\;\forall\;x_0\;\in\mathbb{R}^n\setminus\{0\}$. By substituting $\dot{s}(t,x_0)=As(t,x_0)$ in \eqref{eqn:Ahat} we get the minimum cost $J^*(\hat{A},\;x_0)=0$.
\end{remark}

\subsubsection{Continuity} 
Next, we ask whether one can obtain the local linearization of a nonlinear system around the equilibrium point using \eqref{eqn:Ahat}. The following result shows that this is indeed possible as $\hat{A}(x_0)$ is continuous in its argument $x_0$ at $x_0=0$. The caveat here is that the persistency of excitation condition ($\Gamma_2\succ0$) does not hold at the origin, hence to define continuity of $\hat{A}(x_0)$ at the origin, we need a continuous extension of $\hat{A}(x_0)$, such that $\hat{A}(0)\triangleq A$.  
\begin{theorem}\label{thm:continuity}
   $\lim\limits_{x_0\to0}\hat{A}(x_0)=A$, where
   \begin{align}
    \hat{A}(x_0)= \begin{cases}
    \Gamma_1\Gamma_2^{-1},& x_0\in M\setminus\{ 0\}\\
    A,              &  x_0=0.
\end{cases}
   \end{align}
\end{theorem}
\begin{proof}
    By substituting for $\dot{s}(t,x_0)$ from the dynamics \eqref{eqn:autosystem_lin_ho} in $\Gamma_1$,  the difference $\|\hat{A}(x_0)-A\|$ can be written as $ \|\Gamma_1\Gamma_2^{-1}-A\|=\|\Gamma_3\Gamma_2^{-1}\|$
    where $\Gamma_3\triangleq\int_0^\infty(g(s(t,x_0))s^\top(t,x_0))\mathrm{d}t$.
    Since $g(x)=o(x)$, $\forall\;\epsilon_1>0,\;\exists\;\delta_1>0\;:\;\|x\|<\delta_1\Rightarrow \|g(x)\|<\epsilon_1\|x\|$. Thus, $\|\Gamma_3\Gamma_2^{-1}\|\leq\epsilon_1\Biggr(\int_0^\infty\|s(t,x_0)\|^2\mathrm{d}t\Biggr)\|\Gamma_2^{-1}\|.$
    Next, due to the exponential stability of the origin, $\exists\;\delta_2>0,\;c>0,\;\lambda>0$ such that $\forall\;\|x_0\|<\delta_2,\;\|s(t,x_0)\|\leq c e^{-\lambda t}\|x_0\|$. Hence, $\|\Gamma_3\Gamma_2^{-1}\|\leq\epsilon_1\|x_0\|^2\|\Gamma_2^{-1}\|$.
    Now, one can substitute $s(t,x_0)=e^{At}x_0+d(t,x_0)$, where $d(t,x_0)=o(x_0)$, in the expression for $\Gamma_2$. Next, due to the assumptions that the trajectories are persistently exciting, and the Hurwitzness of $A$, one can prove that $\exists\;k>0,\delta_3>0$ such that $\|\Gamma_2^{-1}\|\leq \dfrac{k}{\|x_0\|^2}, \forall\;\|x_0\|<\delta_3$. We omit this step due to the page limit. Thus we obtain, $ \|\Gamma_3\Gamma_2^{-1}\|\leq\epsilon_1k.$
    Now choosing $\epsilon_1=\epsilon/k$, we arrive at the conclusion that $\forall\;\epsilon>0,\;\exists\;\delta=\mathrm{min}\{\delta_1,\;\delta_2,\;\delta_3\}$ such that $\|x_0\|\Rightarrow\|\hat{A}(x_0)-A\|<\epsilon$.
\end{proof}
\begin{remark}[Finite Length Trajectories]\label{rem: Finitetraj}
    Consider a finite length trajectory, $s(t,x_0):[0,T]\to\mathbb{R}^n$ which is persistently exciting, i.e, $\int_0^Ts(t,x_0)s(t,x_0)^\top \mathrm{d}t\succ0$. Then,  the minimizer of the cost $J(\bar{A},\;x_0)=\int_0^T\|f(s(t,x_0))-\bar{A}s(t,x_0)\|^2\;\mathrm{d}t$, can be obtained by replacing the integration limits in \eqref{eqn:Ahat} with $0$ to $T$.  Thus, one can collect a finite number of system trajectory samples close to the equilibrium to design linear control laws as in \cite{formulabased}, \cite{pasqdatadrivennetworks}. Moreover, the matrix $\hat{A}(x_0)$ as defined in \eqref{eqn:Ahat} can be obtained as $\dot{X}X^\dagger$, by approximating the integrals as Riemann sums where $\dot{X}=\begin{bmatrix}
        \dot{x}(0) & \dot{x}(\delta T)&\cdots&\dot{x}(N\delta T)
    \end{bmatrix}$ and $X=\begin{bmatrix}
        x(0) & x(\delta T)&\cdots&x(N\delta T)
    \end{bmatrix}$, $\delta T$ is the sampling time and $N$ is the number of samples. 
\end{remark}

\subsubsection{Stability of the Minimizer $\hat{A}(x_0)$} 
Now, we prove an important result, that $\hat{A}(x_0)$ is a stable matrix when $x_0\in M$.

\begin{theorem} \label{thm:stability}
    The zero solution $y(t)\equiv0$ of the system $\dot{y}(t)=\hat{A}(x_0)y(t)$ is Lyapunov stable. Furthermore, if the pair $(\hat{A}^\top(x_0),\;x_0^\top)$ is observable, then $\hat{A}(x_0)$ is Hurwitz. 
\end{theorem}
\begin{proof}
    Consider the positive definite matrix $\Gamma_2$, then $\hat{A}(x_0)\Gamma_2+\Gamma_2\hat{A}^\top(x_0) =\Gamma_1+\Gamma_1^\top$.
    Now, the symmetric part $\Gamma_{s} = \Gamma_{1} + \Gamma_{1}^{\top}$ simplifies to
    \begin{align}
        \Gamma_{s}&=\int_0^\infty\Biggr( \dot{s}(t,x_0)s^\top(t,x_0)+s(t,x_0)\dot{s}^\top(t,x_0)\Biggr)\mathrm{d}t\nonumber\\
        &=-x_0x_0^\top \preceq 0.\label{eqn:Lyap_x0} 
    \end{align}
     Inequality \eqref{eqn:Lyap_x0} is obtained using the fundamental theorem of calculus, $\lim\limits_{t\to\infty}s(t,x_0)=0$, and $s(0,\;x_0)=x_0$ as $x_0\in M$. Now, by Theorem 3.18 in \cite{Haddad}, the zero solution $y(t)\equiv0$ is Lyapunov stable. 
    To show that $\hat{A}(x_0)$ is Hurwitz, we invoke Theorem 3.17 in \cite{Haddad}, alongside the assumption that the pair $(\hat{A}^\top(x_0),\;x_0^\top)$ is observable. This concludes the proof.
\end{proof}

\section{Estimation of the Region of Attraction}
\label{sec:results_ROA}
\subsection{Characterization of the Region of Attraction}
To characterize the ROA, we define the \textit{residual energy} as $E:\;\mathbb{R}^{n}\to\mathbb{R}_{\geq0}\cup\{\infty\}$,
\begin{align*}
    E(x_0)\triangleq\lim_{\tau\to\infty}\frac{1}{2}\|\Gamma_1(x_0,\tau)\Gamma_2(x_0,\tau)^{-1}-A\|^2_F
\end{align*}
where $\Gamma_1(x_0,\tau)\triangleq\int_0^\tau\dot{s}(t,x_0)s^\top(t,x_0)\mathrm{d}t$ and $\Gamma_2(x_0,\tau)\triangleq\int_0^\tau s(t,x_0)s^\top(t,x_0)\mathrm{d}t$. 

\begin{lemma}\label{lemma:ROA}
    Consider the dynamics \eqref{eqn:autosystem_lin_ho}, then $0\leq E(x_0)<\infty,\;\forall \; x_0\in M\setminus\{0\}$. Furthermore, if
    \begin{align}
        \lim_{\tau\to\infty}\frac{|\int_{0}^{\tau}g^{\top}(s(t,x_{0}))s(t,x_{0})\mathrm{d}t|}{|\int_{0}^{\tau}s^{\top}(t,x_{0})s(t,x_{0})\mathrm{d}t|} \to \infty \label{eqn: condition on growth}
    \end{align}
    $\forall\; x_0\in \operatorname{int}(M^{c})$ then $E(x_0)=\infty$ when $x_0\in \operatorname{int}(M^{c})$.   
\end{lemma}
\begin{proof}
    First, we show that the residual energy is finite for a given $x_0\in M\setminus\{0\}$. By \eqref{eqn:Ahat},
    \begin{flalign}
        \|\hat{A}(x_0)\|_{F}&=\|\Gamma_1\Gamma_2^{-1}\|_{F}\leq\|\Gamma_1\|_{F}\|\Gamma_2^{-1}\|_{F}\nonumber\\
        &\leq \zeta \Biggr\|\int_0^\infty\dot{s}(t,x_0)s^\top(t,x_0)\mathrm{d}t\Biggr\|_{F}\nonumber\\
        &\leq \zeta \int_0^\infty\|\dot{s}(t,x_0)s^\top(t,x_0)\|_{F}\;\mathrm{d}t\nonumber\\
        &\leq \zeta \Biggr(I_{T}+\int_T^\infty L\|s(t,x_0)\|^2\;\mathrm{d}t\Biggr)\nonumber\\
        &\leq \zeta \Biggr(I_{T} + \int_T^\infty Lce^{-2\mu t}\|x_0\|^2\mathrm{d}t\Biggr)\hspace{-0.1in}\nonumber\\
        &=\zeta(I_{T}+Lc\|x_0\|^2)/2\mu.\nonumber
    \end{flalign}
    where $\displaystyle \zeta=\frac{1}{\lambda_{\mathrm{min}}(\Gamma_2)}$and $I_{T} \displaystyle = \int_{0}^{T}\|\dot{s}(t,x_{0})\|\|s(t,x_{0})\|\mathrm{d}t < \infty$. The third and fourth inequalities are due to the local Lipschitzness of $f(x)$ and local exponential stability, respectively. Now by triangle inequality, $E(x_0)\leq(\|\hat{A}(x_{0})\|_{F}+\|A\|_{F})^{2}$. Thus for a given $x_0\in M\setminus\{0\}$, the residual energy is finite. Next, we analyze the residual energy for $x_0\in\mathrm{int}(M^c)$. Let $x_0\in\operatorname{int}(M^{c})$. Then, define $\hat{A}(x_{0},\tau) = \Gamma_{1}(x_{0},\tau)\Gamma_{2}(x_{0},\tau)^{-1}$, we then have from the definition of $\hat{A}$, equation \eqref{eqn:Ahat}, we get $A\Gamma_2(x_{0},\tau)+\Gamma_3(x_{0},\tau)=\hat{A}(x_0,\tau)\Gamma_2(x_{0},\tau)$, where $\displaystyle \Gamma_3(x_{0},\tau)=\int_0^{\tau} g(s(t,x_0))s^\top(t,x_0)\mathrm{d}t$.
    Re-arranging the terms and using the triangle inequality, $ \|\hat{A}(x_0,\tau)-A\|_{F}\geq\dfrac{\|\Gamma_3(x_{0},\tau)\|_{F}}{\|\Gamma_2(x_{0},\tau)\|_{F}}.$
    The numerator can be written as, 
    \begin{flalign*}
        \hspace{-0.05in}&\|\Gamma_3(x_{0},\tau)\|_{F}=\Biggr\|\int_0^\tau g(s)s^\top \mathrm{d}t\Biggr\|_{F} \hspace{-0.1in}=\sqrt{\sum_{i,j=1}^n\Biggr|\int_0^\tau g_i(s)s_j\;\mathrm{d}t\Biggr|^2}\nonumber\\
         &\geq\sqrt{\frac{1}{n}\Biggr(\sum_{i=1}^n\Biggr|\int_0^\tau g_i(s)s_i\;\mathrm{d}t\Biggr|\Biggr)^2}\geq\frac{1}{\sqrt{n}}\Biggr|\int_0^\tau g^\top(s)s\;\mathrm{d}t\Biggr|.
    \end{flalign*}
    Similarly, the denominator, $\|\Gamma_2(x_{0},\tau)\|_{F}$ is $\|\int_0^\tau ss^\top \mathrm{d}t\|_{F}\leq\int_0^\tau\|ss^\top\|_{F}\;\mathrm{d}t=\int_0^\tau s^\top s\;\mathrm{d}t.$ Then $E(x_{0})=\infty,\;\forall\;x_{0}\in\operatorname{int}(M^{c})$, follows from \eqref{eqn: condition on growth}. 
\end{proof}
\begin{remark}
    Condition \eqref{eqn: condition on growth} can be thought of as defining how the system behaves outside the region of attraction. Specifically, it defines the rate at which the nonlinearities $g(x)$ grow to infinity outside the region of attraction (or at the boundary), along the system trajectories $s(t,x_0)$. This is similar to \cite[Assumption 5]{formulabased}, where the authors assume a (\textit{signal-to-noise ratio}) bound on the nonlinearities to define the behavior of trajectories starting close to the origin. 
\end{remark}

 Lemma \ref{lemma:ROA} enables us to characterize the ROA as a subset of the set of all initial conditions with finite residual energy, i.e., $M\subseteq\{x_0\;\in\mathbb{R}^n\;:\; E(x_0)<\infty\}$. Moreover, if the higher-order terms in $f(x)$ satisfy \eqref{eqn: condition on growth}, the residual energy increases to infinity as $x_0\to\partial M$, much like \textit{barrier functions} discussed in \cite{CBF}. Using these properties, along with continuity of $\hat{A}(x_0)$ as shown in Theorem \ref{thm:continuity}, we propose an algorithm (based on \textit{geometric flows}) to estimate the ROA of systems of the form \eqref{eqn:autosystem} when $f(x)$ is unknown. The algorithm proposed is data-driven and only requires that the trajectories of the system are persistently exciting.
 
\subsection{Algorithm to Estimate the ROA}
By Lemma \ref{lemma:ROA}, $0\leq E(x_0)<\infty\;\forall\;x_0\in M\setminus\{0\}$ and $E(x_0)\to\infty\;\forall\;x_0\to\partial M$. We first map the range of $E(x_0)$ to the interval $[0,1]$, such that $0\leq E(x_0)<1,\;x_0\in M\setminus\{0\}$ and $E(x_0)=1,\;x_0\to\partial M$. Various functions like the hyperbolic tangent, modified exponential and sigmoid functions, etc, achieve this. In this article, we use the hyperbolic tangent $\mathrm{tanh}(E(x_0))$ to map the range of $E(x_0)$ to $[0,1]$. Thus, $\mathrm{tanh}(E(x_0))\in(0,1)\;\forall\;x_0\in M\setminus\{0\}$ and $\mathrm{tanh}(E(x_0))=1\;\forall\;x_0\in \partial M$.

The algorithm converges to the best (exact) estimate of the boundary of the ROA ($\partial M$) by progressively improving an initial estimate. The boundary of the ROA is characterized as an $n-1$-dimensional embedded sub-manifold of $\mathbb{R}^{n}$ and given by $h(x,t) = 0$, where $x\in\mathbb{R}^{n}$ and $t$ denotes time. The following boundary dynamics are proposed, 
\begin{align*}
    \frac{\partial}{\partial t}h(x,t) = [1-\mathrm{tanh}(E(x))]\frac{\partial }{\partial x}h(x,t),\;h(x,0) = h_{0}(x)
\end{align*}
where  $h_0: \mathbb{R}^{n} \to \mathbb{R}$ is a mapping whose zero level set $\partial M_{0} \triangleq h_{0}^{-1}(0)$ denotes the initial conservative estimate  ($M_{0}$) of the ROA. To numerically simulate these dynamics, we propose the following gradient dynamics for a set of representative points $z\in \mathbb{R}^{n}$ satisfying $h(z, t)= 0$. Consider the dynamics
\begin{align}
    \dot{z}(t)=[1-\mathrm{tanh}(E(z(t)))]n(z(t)),\;z(0)\in\partial M_0\label{eqn:geometricflowdynamics}
\end{align}
where $n(z(t))$ is the outward normal at $z(t)$ to $h(z,t)=0$.

\begin{theorem} \label{thm:algoconverge}
    Let $M\subseteq \mathbb{R}^n$ be the ROA of the origin, under the dynamics \eqref{eqn:autosystem_lin_ho}, and let $z(0)\in\partial M_0\subset M$, then the solution of the dynamics \eqref{eqn:geometricflowdynamics} converges to $z^{*}\in \partial M$.
\end{theorem}

\vspace{0.1 cm}
\begin{proof}
       Since $z(0)\in M$ and for all $ z(t)\in M\setminus\{0\}$,  $0<\mathrm{tanh}(E(z(t)))<1$,  $\lim\limits_{t\to\infty}z(t)=z^*$, where $z^*$ satisfies $[1-\mathrm{tanh}(E(z^*))]n(z^*)=0$. Since $n(z^*)\neq0$ as it is an outward normal vector at $z^*$, $\mathrm{tanh}(E(z^*))=1$. Now, by Lemma \ref{lemma:ROA}, this condition is satisfied for points in $\mathrm{int}(M^c)$. Since the dynamics \eqref{eqn:geometricflowdynamics} is initialized using points within $M$,  $z^*$ must belong to $\partial M$. This concludes the proof. 
\end{proof}

\begin{remark}\label{remark:counterexample}
    Though the class of nonlinear systems that satisfy \eqref{eqn: condition on growth} is large, there can be specific systems, for example, $f(x)=\begin{bmatrix}
        -x_1/((1+x_1^2)^2)+x_2&(-x_1-x_2)/(1+x_1^2)^2
    \end{bmatrix}^\top$, that may not satisfy this condition. For such systems, a conservative estimate of the ROA can be constructed by modifying \eqref{eqn:geometricflowdynamics}, as illustrated in Section \ref{sec: simulations}.
\end{remark}

\begin{remark}
    Notice that the dynamics \eqref{eqn:geometricflowdynamics} depends on the value of $E(\cdot)$, which in turn is dependent on $\hat{A}(z)$ and $A$. $\hat{A}(z)$ can be obtained solely based on the state trajectory measurements. Further, the knowledge of $A$ is not required as a local estimate of $\hat{A}(\cdot)$ can be computed using trajectories close to the origin. This estimate can be arbitrarily close to $A$ due to Theorem \ref{thm:continuity} making the algorithm entirely data-driven.
\end{remark} 

\section{Simulations}
\label{sec: simulations}
In this section, we validate the results derived using simulations and estimate the ROA of various systems with unknown dynamics using trajectory data. To critically analyze the results, we provide three examples below. \textbf{a)} system with a bounded ROA, \textbf{b)} system with an unbounded ROA, and \textbf{c)} system not satisfying the radial un-boundedness condition as in Remark~\ref{remark:counterexample}. We encourage the readers to view the simulation results via the following link: https://youtu.be/CBPvZc59XUw .

Since one cannot simulate for an infinite time, we only compute an approximation of $E(x_0)$ with a finite but long enough trajectory. This is achieved by sampling the data at time intervals of $T=0.1$ with a total of $N=40$ samples and replacing infinite integrals in $\hat{A}(x_0)$ with finite Riemann sums, similar to the procedure mentioned in Remark \ref{rem: Finitetraj}. As for the Jacobian matrix $A$, we collect data close to the origin and approximate $A\approx\hat{A}(x_0),\;\|x_0\|<\epsilon$, where $\epsilon=0.1$. This approximation is valid as shown using Theorem \ref{thm:continuity}.

\subsubsection{Bounded Domain of Attraction: Van der Pol Oscillator in Reverse Time}
Consider the following dynamics,
\begin{equation}\label{eqn:simulation_vanderpol}
    \begin{aligned}
         \dot{x}_1=-x_2,\;\dot{x}_2=x_1-(1-x_1^2)x_2.
    \end{aligned}
\end{equation}

\begin{figure}
        \centering
        \includegraphics[width=0.39\textwidth]{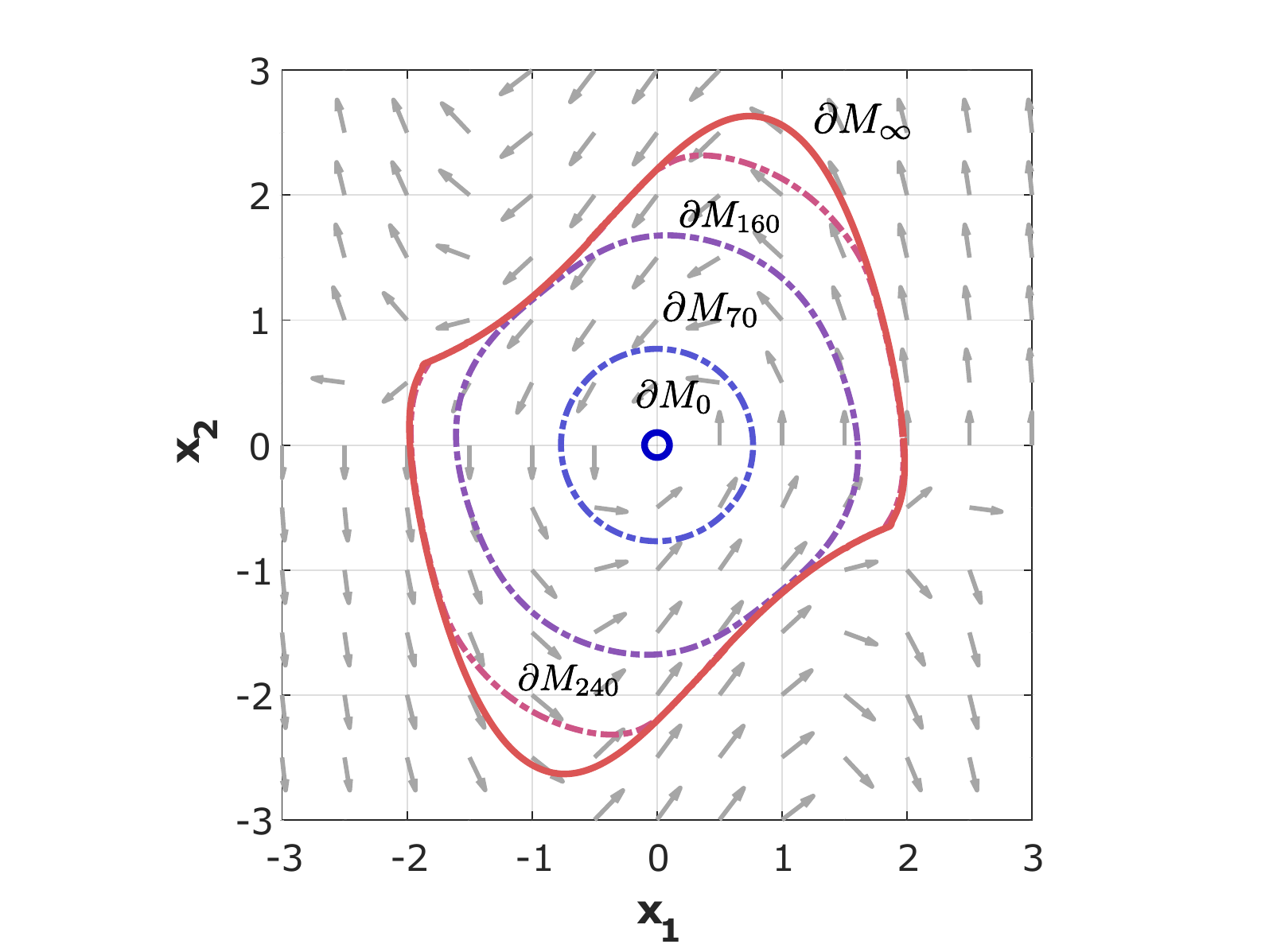}
        \caption{The evolution of the boundary of the ROA is shown, where $\partial M_i$ represents the boundary submanifold at the $i^{\mathrm{th}}$ iteration. The algorithm is initialized with the curve $\partial M_0$ and is shown to asymptotically converge to $\partial M_\infty$, which is the true boundary of the ROA.}
        \label{fig: vander}
    \end{figure}

The origin of \eqref{eqn:simulation_vanderpol} is a locally exponentially stable equilibrium point with a region of attraction bounded by an unstable limit cycle. We sample $D=50$ points from an initial guess of the boundary of the region of attraction $\partial M_0=\{x\in\mathbb{R}^2:\;\|x\|=0.1\}$ (or a circle or radius $0.1$). Then, these 40 sample points are propagated forward under the dynamics \eqref{eqn:geometricflowdynamics}. The evolution of the closed curve $\partial M_t$ is shown in Fig~\ref{fig: vander}. We can observe that the curve $\partial M_t$ converges to the exact boundary of the ROA(the limit cycle).  
\subsubsection{Conservative Estimate of an Unbounded ROA}
\begin{figure}
        \centering
        \includegraphics[scale=0.28]{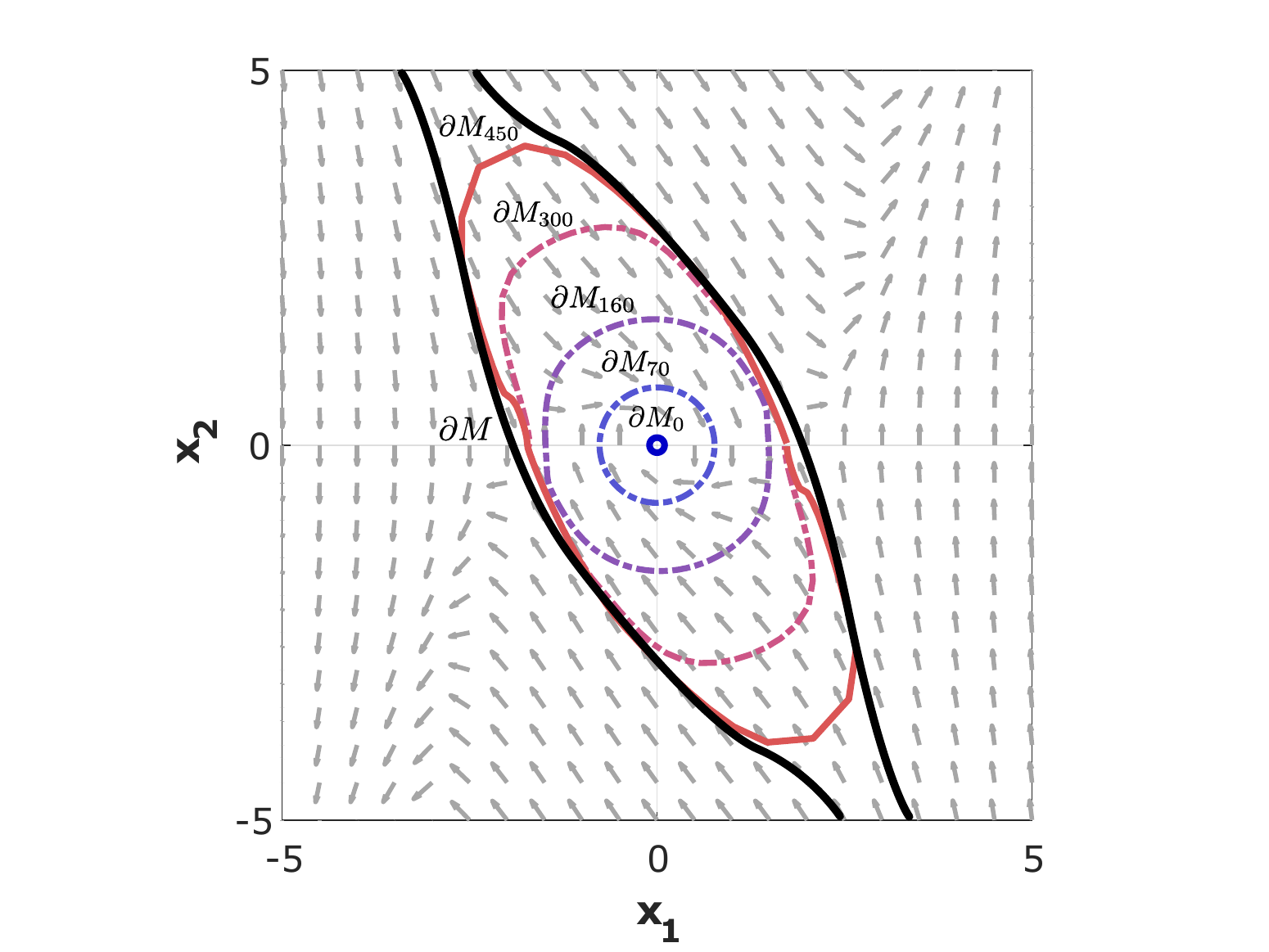}
        \caption{The evolution of $\partial M_i$ for a system with an unbounded ROA. The solid black curves represent the actual ROA. Since the algorithm was terminated in finite time, a conservative estimate $\partial M_{450}$ is obtained.}
        \label{fig:Unbounded}
\end{figure}

Consider the dynamical system given by
\begin{equation}\label{eqn:simulation_unbounded}
    \begin{aligned}
         \dot{x}_1&=x_2,\;\dot{x}_2=-x_1-x_2+\frac{1}{3}x_1^3.
    \end{aligned}
\end{equation}
The domain of attraction for this system is unbounded. Hence, we can only construct a conservative estimate of the ROA in finite time. The actual and the estimated domain of attraction has been shown in Fig~\ref{fig:Unbounded}. Note that the estimate $\partial M_t$ would have converged to the actual boundary of the domain of attraction asymptotically if run for an infinite time. 

Both the systems \eqref{eqn:simulation_vanderpol} and \eqref{eqn:simulation_unbounded} satisfy the condition \eqref{eqn: condition on growth}, but there could be systems as mentioned in Remark \ref{remark:counterexample} where the \textit{growth rate} of $g(x)$, outside the region of attraction, is finite. One can still obtain a conservative estimate of the ROA of such systems by modifying the dynamics as
\begin{align}\label{eqn:modifiedalgo}
    \dot{z}(t)=[\gamma-\mathrm{tanh}(E(z(t)))]n(z(t)),\;z(0)\in \mathcal{B}_{\epsilon}(0)
\end{align}
where $\mathcal{B}_{\epsilon}(0)$  is an $\epsilon$-ball around the origin and $\gamma\in(0,1)$. Choosing a $\gamma\in(0,1)$ will make sure the dynamics \eqref{eqn:modifiedalgo} converges to a point $z^*\in M$. (Note: A smaller $\gamma$ will give a more conservative estimate of the ROA. But choosing a larger $\gamma$ will lead to errors in the estimate of the ROA as one cannot guarantee that the algorithm will converge to points inside the ROA in this case.) Thus, we can construct a conservative estimate of the boundary of the ROA as $\partial M_c=\{z^*\}\subset M$.  We illustrate this using the next example.
\subsubsection{Conservative Estimate of the ROA of a System which does not satisfy \eqref{eqn: condition on growth}}
\begin{figure}
        \centering
        \includegraphics[scale=0.17]{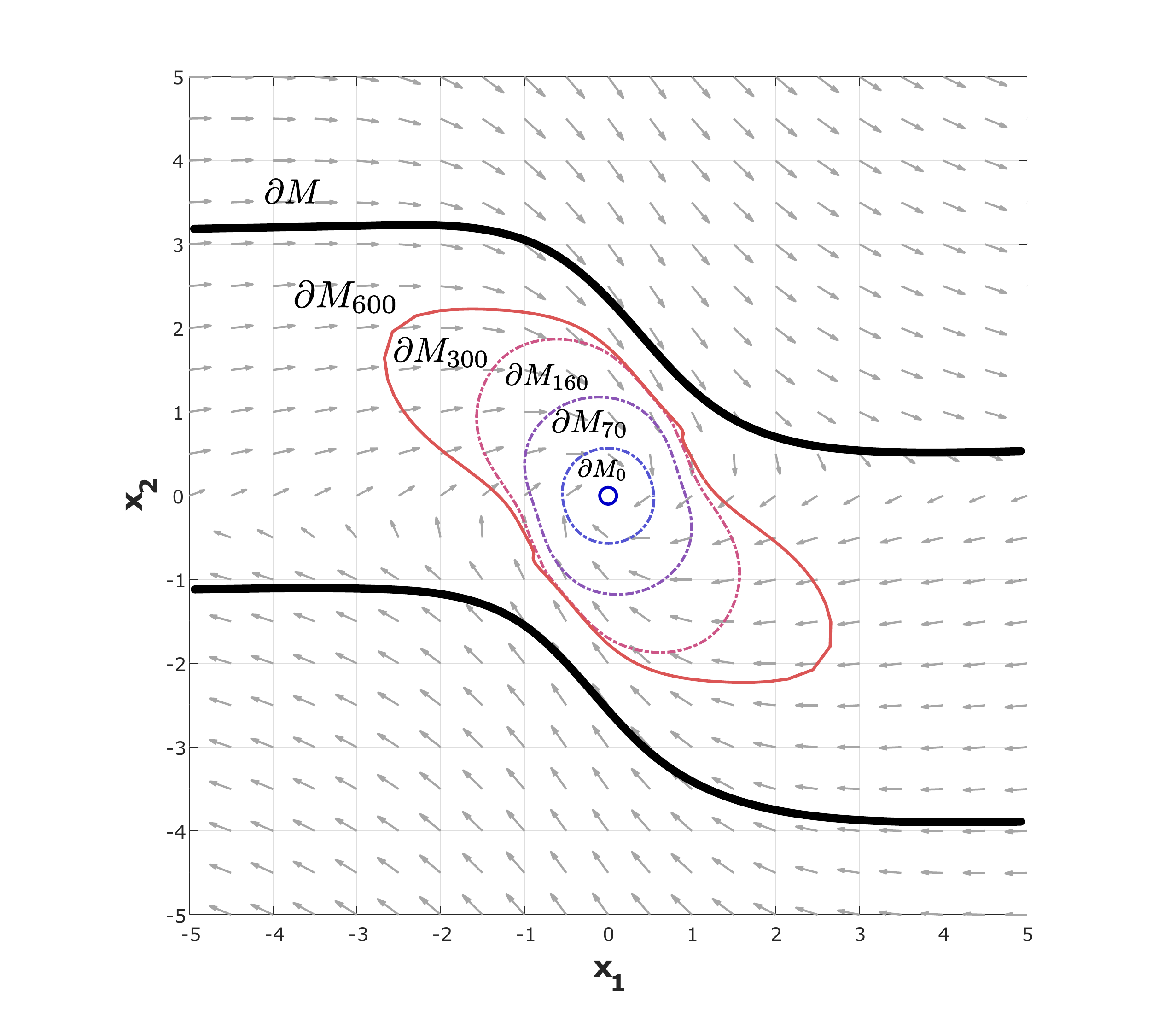}
        \caption{The system has an unbounded ROA, but assumption \eqref{eqn: condition on growth} is not satisfied. Hence, only a conservative estimate can be obtained by setting $\gamma=0.7$. The actual boundary of the ROA is shown as solid black curves.}
        \label{fig:rational}
    \end{figure}

Consider the system dynamics
\begin{equation*}
    \begin{aligned}
        \dot{x}_1=[-x_1/(1+x_1^2)^2]+x_2,\;\dot{x}_2=-(x_1+x_2)/(1+x_1^2)^2.
    \end{aligned}
\end{equation*}
Approximating the system behavior for large values of $\|x(t)\|^2$, and after some simple calculations, one can prove that the limit given in \eqref{eqn: condition on growth} is finite in this case. Hence, to estimate the ROA, we set $\gamma=0.7$ in \eqref{eqn:modifiedalgo} and obtain a conservative estimate of the ROA, as shown in Fig~\ref{fig:rational}. 

\section{Conclusion and Future Work}
\label{sec:conclusion_Future}
The article introduced a few results on analog data-driven control theory. We established that a global minimizer exists for the error function between the system dynamics and a linear estimate. Furthermore, for initial conditions inside the region of attraction, we proved that this global minimizer is Lyapunov stable and, under certain conditions, Hurwitz. We also provide an algorithm to compute the estimate of the boundary of the region of attraction without prior knowledge of the system model. The algorithm improves the initial conservative guess of the boundary of the region of attraction and converges to the best (exact) estimate. In our future work, we aim to expand the results for systems with control inputs. We are also working on extending the results to estimate the region of attraction on manifolds using geometric flows.
\bibliographystyle{IEEEtran}
\bibliography{ref}
\end{document}